\documentclass[12pt,a4paper]{amsart} 
\usepackage{geometry}
\usepackage{pinlabel}         
\usepackage{graphicx}
\usepackage{amssymb, amsthm}
\usepackage{epstopdf}
\def\Z{\mathbb{Z}}

\def\l{\lambda}

\def\<{\langle}
\def\>{\rangle}
\def\fix{{\rm{Fix}}}

\def\ro{R}
\def\e{e}

\def\autn{{\rm{Aut}}$(F_n)$}
\def\sautn{{\rm{SAut}}$(F_n)$}
\def\sln{${\rm{SL}}(n,\mathbb Z)$}
\def\Autn{{\rm{Aut}}(F_n)}

\def\Sautn{{\rm{SAut}}(F_n)}
\def\gln{${\rm{GL}}(n,\mathbb Z)$}
\def\GL{{\rm{GL}}(n,\mathbb Z)}
\def\Z{\mathbb Z}
\def\R{\mathbb R}
\def\bS{\mathbb S}

\def\SN{SN}   
\def\SW{SW_n}
\newcommand{\A}{\mathcal A}
\newcommand{\Oh}{\mathcal O}
\newcommand{\Hcp}{H^c}
\newcommand{\Hc}{H_c}

\newtheorem{theorem}{Theorem}[section]
\newtheorem{thm}[theorem]{Theorem}
\newtheorem{lemma}[theorem]{Lemma}
\newtheorem{prop}[theorem]{Proposition}
\newtheorem{corollary}[theorem]{Corollary}
\newtheorem{definition}[theorem]{Definition}
\newtheorem{remark}[theorem]{Remark}

\begin{document}

\catcode`\@=11
\def\serieslogo@{\relax}
\def\@setcopyright{\relax}
\catcode`\@=12

\def\jump{\vskip 0.3cm}
\title[Actions of
${\rm{Aut}}(F_n)$ on homology manifolds]{Actions of automorphism groups of free
groups on homology spheres and acyclic manifolds}

\author[Bridson]{Martin R.~Bridson}
\address{Martin R.~Bridson\\
Mathematical Institute \\
24-29 St Giles' \\
Oxford OX1 3LB \\
U.K. }
\email{bridson@maths.ox.ac.uk}

\author[Vogtmann]{Karen Vogtmann }
\address{Karen Vogtmann\\
Department of Mathematics\\
Cornell University\\
Ithaca NY 14853 }
\email{vogtmann@math.cornell.edu}
 
\thanks{Bridson's work was supported by an EPSRC Senior Fellowship
and a Royal Society Nuffield Research Merit
Award. Vogtmann is supported NSF grant DMS-0204185. Bridson thanks the
University of Geneva and the EPFL (Lausanne) for their hospitality during the writing
of this paper, and the Swiss National Research Foundation for funding his stay.}

\subjclass{20F65, 20F28, 53C24, 57S25} 


\keywords{automorphism groups of free groups, rigidity,
group actions, Smith theory, generalized manifolds}

\begin{abstract} For $n\ge 3$, let \sautn\  denote the unique subgroup of index two in the
automorphism group of a free group. 
The standard linear action
of \sln\  on $\R^n$ induces  non-trivial actions of \sautn\ 
on $\R^n$ and on
$\mathbb S^{n-1}$.
We prove that 
\sautn\ admits no non-trivial actions by homeomorphisms
on acyclic manifolds or spheres of smaller dimension. Indeed,  
\sautn\ cannot act non-trivially on any generalized
$\Z_2$-homology sphere of dimension less than $n-1$, nor on any
$\Z_2$-acyclic $\Z_2$-homology manifold of dimension less than $n$. It follows
that \sln\ cannot act non-trivially on such spaces either. When  $n$ is even,
we obtain similar results with $\Z_3$ coefficients.  
\end{abstract}

\maketitle

\section{Introduction} 

In geometric group theory one attempts to elucidate the algebraic properties of a group by studying its actions on spaces with good geometric properties.  For irreducible lattices in higher-rank semisimple Lie groups, versions of
Margulis superrigidity  place severe restrictions on the spaces that
are useful for this purpose. 
Our focus in this article is on the rigidity properties of
the group $\Autn$ of automorphisms of a free group, which is not a lattice  but neverthless enjoys many similar properties.

In
\cite{BriVog03:2} we exhibited strong constraints on homomorphisms from \autn\ 
and pointed out that such constraints restrict the way in which  \autn\ 
can act on various spaces. 
We illustrated this
point by showing that if $n\ge 3$ then any  action of \autn\  
on the circle by homeomorphisms must factor through the determinant
homomorphism $\det\colon \Autn\to \Z_2$.  We now show that similar resrictions
apply much more generally, to  actions on higher-dimensional generalized homology spheres over  $\Z_p$
  and to generalized manifolds that are  $\Z_p$-acyclic, for  $p=2,3$.  

For $n\geq 3$ we denote by  \sautn\ the unique subgroup of index two in 
\autn. The action of \autn\ on the abelianization of the free group $F_n$
gives a natural map $\Autn\to \GL$, sending \sautn\  onto \sln.   Thus
the standard linear action
of \sln\  on $\R^n$ induces  non-trivial actions of \sautn\ 
on $\R^n$ and on the sphere 
$\mathbb S^{n-1}$. However, we will
 prove that \sautn\ cannot act non-trivially
on spheres or contractible manifolds of any smaller dimension.   
For linear actions, elementary results in the representation
theory of finite groups can be combined with an understanding of the torsion
in \sautn\ to prove this statement; the real challenge lies with non-linear
actions.

Smooth actions are considerably easier to handle than topological ones.
 Thus we begin by proving, in
section \ref{s:smooth}, that  for $n\ge 3,$ 
\sautn\ cannot act non-trivially by diffeomorphisms on a $\Z_2$-acyclic
smooth manifold of dimension
less than $n$.  The proof we present is deliberately constructed so as to point out the difficulties encountered in the purely topological setting.    
In particular, the proof requires understanding the fixed point sets of involutions.  This
 immediately creates a problem in the topological setting because the fixed point
sets  of involutions are not in general manifolds, but only homology manifolds over $\Z_2$.   
A second difficulty arises because there is no tangent space in the topological setting; in the smooth case the tangent space allows one to use linear algebra to transport 
information about the action near fixed point sets 
to information about the action on the ambient manifold. 

These are
well-known difficulties that lie at the heart of the theory of transformation groups and
much effort has gone into confronting them \cite{Borel}, \cite{BredonBook}.
They are overcome using (local and global) {\em{Smith theory}}, but one
has to accept the necessity of working with generalized manifolds rather than
classical manifolds. (See section \ref{s:sphere} for definitions concerning generalized
manifolds.)

We shall prove the following results by following the architecture of the proof we
give in the smooth setting, combining Smith theory with an
analysis of the torsion  in \sautn\ to overcome the technical problems that arise.

\begin{thm}\label{t:2sphere}
If $n\ge 3$ and $d<n-1$, then any action of \sautn\  
by homeomorphisms on a generalized $d$-sphere over $\Z_2$
 is trivial, and hence \autn\ can act  only via
the determinant map.
\end{thm}

\begin{thm}\label{t:2acyclic}
If $n\ge 3$ and $d<n$, then any action of \sautn\  
by homeomorphisms on a $d$-dimensional
$\Z_2$-acyclic homology manifold over $\Z_2$
 is trivial, and hence \autn\ can act  only via
the determinant map.
\end{thm}

As special cases we obtain the desired minimality result for the standard
linear action of \sautn\ on $\mathbb R^n$ and $\mathbb S^{n-1}$.

\begin{corollary} If $n\ge 3$,   then 
\sautn\ cannot act non-trivially by homeomorphisms on 
any contractible manifold of dimension less than $n$,
or on any sphere
of dimension less than $n-1$.
\end{corollary}

We also note that these theorems  have as immediate corollaries the analogous statements for \sln\ and \gln.

\begin{corollary}\label{c:sl}
If $n\ge 3$ and $d<n$, then \sln\  cannot act non-trivially
by homeomorphisms on any generalized $(d-1)$-sphere over $\Z_2$,
 or on any $d$-dimensional
 homology manifold over $\Z_2$ that is $\Z_2$-acyclic.
 Hence \gln\  can act on such spaces only via
the determinant map.
\end{corollary}

Corollary \ref{c:sl} was conjectured by Parwani \cite{Par06}; see remark
\ref{ParFlaw}. 

In section \ref{s:quots} we describe a subgroup $T\subset{\rm{SAut}}(F_{2m})$  
 isomorphic to $(\Z_3)^m$ that intersects every proper
normal subgroup of ${\rm{SAut}}(F_{2m})$ 
trivially. This provides a stronger degree of rigidity than is offered by the 2-torsion
in \sautn\ and consequently one can deduce the following theorems from Smith theory more
readily than is possible in the case of $\Z_2$ (see section \ref{ss:Zthree}).

\begin{thm}\label{t:3sphere}
If $n> 3$ is even and $d<n-1$, then any action of \sautn\  
by homemorphisms on a generalized $d$-sphere over $\Z_3$  is trivial.
\end{thm}

\begin{thm}\label{t:3acyclic}
If $n> 3$ is even and $d<n$, then any action of \sautn\  
by homemorphisms on a $d$-dimensional
$\Z_3$-acyclic homology manifold over $\Z_3$ is trivial.
\end{thm}

 We expect that our results concerning \sln\ should be true for 
 other lattices in ${\rm{SL}}(n,\mathbb R)$, 
 but our techniques do not apply because we make essential use of the torsion in \sln.
 What happens for subgroups of finite index in \sautn\ is less clear: there are
 subgroups of finite index in \sautn\ that map non-trivially to 
 ${\rm{SL}}(n-1,\mathbb R)$ and hence act non-trivially on $\R^{n-1}$, but
 one does not know if such subgroups can act non-trivially on contractible
 manifolds of dimension less than $n-1$. 
 \smallskip
 
In a brief final section we explain how our results
concerning  torsion in \autn, together with the application
of Smith theory in \cite{MannSu}, imply the following result. 

\begin{thm}\label{noAct} Let $p$ be a prime and let $M$ be a compact $d$-dimensional 
 homology
manifold  over $\Z_p$.   There
exists an integer $\eta(p,d,B)$, depending only on $p,d$ and the sum $B$ of the   mod $p$ betti numbers of $M$,
so that  \autn\ cannot act
non-trivially by homeomorphisms on $M$ if $n>\eta(p,d,B)$.
\end{thm}
  
\smallskip

{\bf Acknowledgements}. We would like to thank the colleagues 
who helped us struggle with the technicalities of generalized manifolds and Smith theory over the past year, including in particular Mladen Bestvina, Mike Davis, Ian Hambleton and Shmuel Weinberger. We also thank Linus Kramer and Olga Varghese for their comments concerning Theorem 4.5.

\section{Smooth actions}\label{s:smooth}

In this section we indicate how to prove Theorem \ref{t:2acyclic} for smooth actions. Our intent here is to explain the structure of the proof of our general results without the technical difficulties that occur in the topological setting.

\begin{theorem}\label{t:smooth}
Let $X$ be a $k$-dimensional differentiable manifold that is $\Z_2$-acyclic (i.e. has the $\Z_2$-homology of a
point). If $n\ge 3$ and $k<n$ then any action of 
\sautn\  by diffeomorphisms on $X$ is trivial.
\end{theorem}

\newcommand{\eij}{\varepsilon_{ij}}
\newcommand{\eab}{\varepsilon_{12}}
\newcommand{\ebc}{\varepsilon_{23}}
\newcommand{\eac}{\varepsilon_{13}}
\newcommand{\eyz}{\varepsilon_{45}}

\begin{proof}  The proof proceeds by induction on $n$.   We omit the cases $n 
\leq 4$, where ad hoc arguments apply (cf. subsection 4.5). Suppose, then, that $n 
\geq 5$, fix a basis 
$a_1,\ldots, a_n$ for $F_n$ and consider the involutions $\eij$ of $F_n$ defined as follows:  
 \[
 \eij\colon\begin{cases} 
            a_i\mapsto a_i^{-1}\cr
            a_j\mapsto a_j^{-1}\cr
            a_k\mapsto a_k& k\neq i,j
            \end{cases}
\]
These involutions are all conjugate in \sautn, and the quotient of \sautn\ by the normal closure of any $\eij$ is ${\rm{SL}}(n,\Z_2)$, which is a simple group (cf. Proposition~\ref{p:quots}).  Thus to prove that an action of \sautn\  is trivial it suffices to  show first that some $\eij$ acts trivially, so that the action factors through ${\rm{SL}}(n,\Z_2)$,  and then that some non-trivial element of ${\rm{SL}}(n,\Z_2)$ acts  trivially. 

Since $X$ is $\Z_2$-acyclic it must be orientable, and since \sautn\ is perfect it must act by orientation-preserving diffeomorphisms.  Therefore either the action of $\eab$ is trivial or the fixed point set $\mathcal F_{12}$ of $\eab$ is a smooth submanifold of codimension at least 2, and Smith theory [22] tells us that this fixed point set will itself be  $\Z_2$ -acyclic.  

The centralizer of $\eab$ contains an obvious copy of ${\rm{SAut}}(F_{n-2})$,
corresponding to the sub-basis $a_3,\dots,a_n$, and by induction this must act
trivially on $\mathcal F_{12}$.  In particular, the automorphism $\eyz$ acts trivially, so its fixed point set $\mathcal F_{45}$ contains $\mathcal F_{12}$.  But $\eab$ and $\eyz$ are conjugate, so in fact $\mathcal F_{12}=\mathcal F_{45}$, i.e. we have two commuting involutions with the same (non-empty) fixed point set.  On the tangent space at a common fixed point these induce commuting linear involutions of $\mathbb R^k$ with the same fixed vectors, which must be identical by basic linear algebra.  But the action of a finite group on a connected
smooth manifold is determined by its action on the tangent space of a fixed point, so the actions themselves must be identical.   Thus the product $\eab\eyz$ acts trivially.  A similar argument shows that $\ebc\eyz$ acts trivially, and we conclude that the product $\eab\eyz\ebc\eyz=\eac$ acts trivially.

Now look at the induced action of  ${\rm{SL}}(n,\Z_2)$ on $X$, and consider the elementary matrices $E_{1j}$. These generate a subgroup isomorphic to $\Z_2^{n-1}$, and we claim that any such group acting by orientation-preserving homeomorphisms on $X$ must contain an element which acts trivially.  To see this, choose an element of $\Z_2^{n-1}$ whose fixed point set $\mathcal F$ has the largest dimension.  By induction (starting with the trivial case $n=3$), some other element of the group must act trivially on $\mathcal F$, and one thus obtains two commuting involutions that have the same fixed point set, as in the previous paragraph.  As before, the involutions must be the same and the product acts trivially.  
\end{proof}
 
\section{Concerning the quotients of \autn}\label{s:quots}

\subsection{Notation}

Fix a generating set $\{a_1,\ldots,a_n\}$ for $F_n$. The right and left Nielsen automorphisms $\rho_{ij}$ and $\lambda_{ij}$ are defined by
\[
{\rho_{ij}\colon\begin{cases}  a_i\mapsto a_ia_j&\cr
            a_k\mapsto a_k& k\neq i\end{cases}}
            \qquad\qquad{\lambda_{ij}\colon\begin{cases}  a_i\mapsto a_ja_i\cr
            a_k\mapsto a_k& k\neq i\end{cases}}
\]
We denote by $e_i$ the automorphism which inverts the generator  $a_i$.  Elements of the subgroup $\Sigma_n$ of automorphisms which permute the generators $a_i$ will be denoted using standard cycle notation; for example $(ij)$ is the automorphism interchanging $a_i$ and $a_j$.
\[
{e_i\colon\begin{cases}  a_i\mapsto a_i^{-1}&\cr
            a_k\mapsto a_k & k\neq i
            \end{cases}}
\qquad\qquad{(ij)\colon\begin{cases}  a_i\mapsto a_j\cr
            a_j\mapsto a_i\cr
            a_k\mapsto a_k& k\neq i,j\end{cases}}
\]
$W_n$ is the subgroup of \autn\ generated by $\Sigma_n$ and the inversions $e_i$, and $\SW$ is the intersection of $W_n$ with \sautn.  The subgroup of $W_n$ generated by the $e_i$ is a normal subgroup $N\cong (\Z_2)^n$, and $W_n$ decomposes as the semidirect product $N\rtimes \Sigma_n$.  The intersection of $N$ with \sautn\ is denoted $\SN$.   Note that the central element $\Delta=e_1e_2\cdots e_n$ of $W_n$ is in $\SN$ if and only if $n$ is even.

Although it seems awkward at first glance,
it is convenient to work with the right action of \autn\ on $F_n$: so
$\alpha\beta$ acts as $\alpha$ followed by $\beta$. An advantage of
this is the neatness of the formula $[\lambda_{ij},\lambda_{jk}]=\lambda_{ik}$,
where our commutator convention is $[a,b]=aba^{-1}b^{-1}$.

\subsection{How kernels can intersect $\SW$}

The following variation on Proposition 9 of \cite{BriVog03:2} will
be useful here. 

\begin{prop}\label{p:quots}
Suppose $n\ge 3$ and let $\phi$ be a  homomorphism from \sautn\ to
a group $G$.  If $\phi|_{\SW}$ has non-trivial kernel $K$, then one of the following holds:

1. $n$ is even, $K=\langle \Delta\rangle$ and $\phi$
factors through ${\rm{PSL}}(n,\Z)$, 

2. $K=\SN$ and  the image of $\phi$ is isomorphic to ${\rm{SL}}(n,\Z_2)$,  or

3. $\phi$ is the trivial map.

\end{prop}

\begin{proof}
   In \autn\ one has the semidirect product
decomposition $W_n=N\rtimes \Sigma_n$ and accordingly we write
elements of $SW_n$ as $\alpha\sigma$, with $\alpha=e_1^{\epsilon_1}e_2^{\epsilon_2}\ldots e_n^{\epsilon_n}\in N$ and $\sigma \in S_n$. (Note that it
may be that neither $\alpha$ nor $\sigma$ is  itself in \sautn.)

Using exponential notation to denote conjugation, we have
\begin{equation}\label{conjugation}
\l_{ij}^{\alpha}=
\begin{cases}
\l_{ij} & \text{if $\epsilon_i=\epsilon_j=0$,}
\\
\rho_{ij} & \text{if $\epsilon_i=\epsilon_j=1$,}
\\
\l_{ij}^{-1} & \text{if $\epsilon_i=0$ and $\epsilon_j=1$,}
\\
\rho_{ij}^{-1} & \text{if $\epsilon_i=1$ and $\epsilon_j=0$.}
\end{cases}
\end{equation}
Also, for $\theta\in\{\l,\rho\}$, we have $\theta_{ij}^{\sigma} = \theta_{\sigma(i)\sigma(j)}$.
Hence $\l_{ij}^{\alpha\sigma} = \theta_{\sigma(i)\sigma(j)}^{\pm 1}$ for some
for $\theta\in\{\l,\rho\}$.

 If  $K$ contains the center $\langle\Delta\rangle$ of $W_n$ then $n$ must be even and the relations
$\Delta \lambda_{ij} \Delta = \rho_{ij}$ imply that the map $\phi$
factors through \sln, since  by \cite{Ger84} adding the relations $\lambda_{ij}=\rho_{ij}$ to a presentation for \sautn\ gives a presentation for \sln. Since $\Delta$ maps to the center of \sln, the map in fact factors through ${\rm{PSL}}(n,\Z)$.  

If $K$ contains an element $\alpha\in\SN$ which is not central in $W_n$, then we can write $\alpha=e_1^{\epsilon_1}e_2^{\epsilon_2}\ldots e_n^{\epsilon_n}$ with $\sum \epsilon_i$  even and some $\epsilon_k=0$.  Given any indices $i$ and $j$ we can conjugate $\alpha$ by an element of the alternating group $A_n\leq \SW$ to obtain elements in the kernel of $\phi$ with any desired values of $\epsilon_i,\epsilon_j\in\{0,1\}$.  Conjugating $\lambda_{ij}$ by these elements, we see from (1) that $\lambda_{ij},\rho_{ij},\lambda_{ij}^{-1}$ and $\rho_{ij}^{-1}$ all have the same image under $\phi$. This implies not only that $\phi$ factors  through ${\rm{SL}}(n,\Z)$, but also that the images of all Nielsen automorphisms have order 2, and so $\phi$ factors through  ${\rm{SL}}(n,\Z_2)$. The image of $\SN$ is trivial under this map, i.e. $K\supseteq \SN$.  Since ${\rm{SL}}(n,\Z_2)$ is simple, the image of $\phi$ is either trivial or isomorphic to  $\rm{SL}(n,\Z_2)$.

Finally, suppose that $K$ contains an element $\alpha\sigma$ which is not in $\SN$, i.e. $\sigma\neq 1$. 
If $\sigma$ is not an involution, then for some $i,j,k$
with $i\neq k$ we have $\sigma(i)=j$ and $\sigma(j)=k$,
hence $\l_{ij}^\sigma = \theta_{jk}$ with $\theta\in\{\lambda,\rho\}$.
By combining the relations $[\l_{ij},\l_{jk}]=\l_{ik}$ and $[\theta_{ij}^\pm,\l_{ij}]=1$
with the fact that $\phi(x^{y})=\phi(x)$ for all
$x\in\Sautn$ and $y\in K$, we deduce:
\begin{align*}
\phi(\l_{ik}) &= [\phi(\l_{ij}),\phi(\l_{jk})] =
[\phi(\l_{ij}^{\alpha\sigma}),\phi(\l_{jk})]\\ & = [\phi(\theta_{jk}^{\pm 1}),\phi(\l_{jk})] =
\phi([\theta_{jk}^{\pm 1},\l_{jk}] )=1.
\end{align*}
Since all Nielsen automorphisms are conjugate in \sautn\ and they together
generate \sautn, we conclude that $\phi$ is trivial (and $K=\SW$).

Finally,
if $\sigma$ is an involution interchanging $j\neq k$, then a similar
calculation produces the conclusion that
$$
\phi(\l_{ik}) =
[\phi(\l_{ij}),\phi(\l_{jk}^{\alpha\sigma})] =
\phi([\l_{ij}, \theta_{kj}^{\pm 1}]) = 1
$$
so that $\phi$ is again trivial.  
\end{proof}

\subsection{All non-trivial quotients of ${\rm{SAut}}(F_{2m})$ contain $(\Z_3)^m$}

In this subsection we are only interested in free groups of even rank. It
is convenient to switch notation: if
$n=2m$ we fix a basis $\{a_1,b_1,\dots,a_m,b_m\}$ for $F_n$; we write
$\lambda_{a_ib_i}$ and $\rho_{a_ib_i}$ for the Nielsen
transformations that send $a_i$ to $b_ia_i$ and $a_ib_i$, respectively;
we write $(a_i\ b_i)$ for the automorphism that interchanges $a_i$
and $b_i$, fixing the other basis elements; we
write $e_{a_1}$ instead of $e_1$, and so on.

\medskip

Let $T$ be the subgroup of \sautn\ generated by 
$\{\ro_i\mid i=1,\dots,m\}$ where  
\[
{\ro_{i}\colon\begin{cases}  a_i\mapsto b_i^{-1} &\cr
            b_i\mapsto b_i^{-1}a_i& \cr
            a_j\mapsto a_j& j\neq i\cr
            b_j\mapsto  b_j& j\neq i.\end{cases}} 
\]

\begin{lemma}\label{whatTis} $T\cong (\Z_3)^m$.
\end{lemma}

\begin{proof} One can verfiy this by direct calculation but the nature
of $T$ is most naturally described in terms of the  labelled  graph $\mathcal T_m$
depicted in Figure~\ref{Tgraph}.

\begin{figure}[ht!]
\labellist
\small\hair 2pt  
\pinlabel {$v_0$} at 105 -7
\pinlabel {$a_1$} at 42 28
\pinlabel {$b_1$} at 62 49
\pinlabel {$v_1$} [r] at 30 65
\pinlabel {$a_2$} at 60 62
\pinlabel {$b_2$} at 87 75
\pinlabel {$v_2$} [r] at 72 95
\pinlabel {$a_3$} at 95 65
\pinlabel {$b_3$} at 125 62
\pinlabel {$v_3$} [r] at 120 102
\pinlabel {$a_m$} at 142 46
\pinlabel {$b_m$} at 164 21
\pinlabel {$v_m$} [r] at 195 60

\endlabellist
\centering
\includegraphics[scale=1.0]{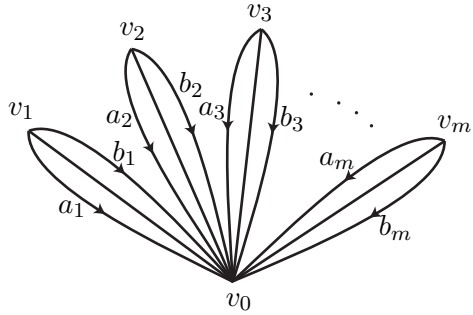}
\caption{Graph realizing the subgroup $T$} 
\label{Tgraph}
\end{figure}

$\mathcal T_m$ has   $m+1$ vertices $v_0,v_1,\dots,v_m$ and
 $3$ edges joining $v_0$ to each of the other vertices. A maximal
tree is obtained by choosing an (unlabelled) edge joining $v_0$ to each of the other vertices.
For each $i$, the remaining two edges incident at $v_i$ are oriented towards
$v_0$ and labelled $a_i$ and $b_i$ .  

This labelling  identifies $\pi_1(\mathcal T_m,v_0)$ with
$F_{2m}=F(a_1,b_1,\dots,a_m,b_m)$ and defines an injective homomorphism
$\psi: {\rm{Sym}}(\mathcal T_m,v_0)\to \Autn$ whose image contains $T$.
Indeed $\ro_i$ is the image under $\psi$ of the symmetry of order 3 that
cyclically permutes the edges joining $v_i$ to $v_0$, sending the edge labelled
$a_i$ to that labelled $b_i$ and sending the edge labelled $b_i$ to the
unlabelled edge.
\end{proof}

A routine calculation yields:

\begin{lemma}\label{3lem} For  $i=1,\dots,m$, let $\beta_i\in\Sautn$  be the automorphism that sends
$a_i$ to $a_i^{-1}$ and
$b_i$ to $a_i^{-1}b_i^{-1}a_i$ while fixing the other basis elements.
\begin{enumerate}
\item  $R_i\e_{a_i}\e_{b_i}R_i^{-1}=\beta_i$.
\item $[R_j,\e_{a_i}]=[R_j,\e_{b_i}]=1$ if $j\neq i$.
\item  $R_i\e_{b_i} R_i^{-1}\e_{a_i}=\lambda_{b_ia_i}^2$.
\item  $R_i^{-1}\e_{a_i} R_i\e_{b_i}=\lambda_{a_ib_i}^2$.
\end{enumerate}
\end{lemma}

\begin{prop}\label{p:3quots} For $m\ge 2$ and any group $G$, let  $\phi\colon{\rm{SAut}}(F_{2m})\to G$ be a homomorphism. If $\phi|_T$ is not injective, then $\phi$ is trivial. \end{prop} 

\begin{proof}  
Let $t\in T$ be a non-trivial element of the kernel of $\phi$.  Replacing $t$ by $t^{-1}$ if necessary, we may write 
$t=\ro_i  u$ where $u$ is a word in the $\ro_j$ with $j\neq i$.
Since each $R_j$ commutes with $e_{a_i}$ and $e_{b_i}$,  we have $t e_{a_i} e_{b_i} t^{-1} = R_i\e_{a_i}\e_{b_i}R_i^{-1}=\beta_i.$  
Since $\phi(t)=1$, applying $\phi$ to this equation gives $\phi(e_{a_i} e_{b_i})=\phi(\beta)$.

We now note that $\e_{a_i}\e_{b_i}$ conjugates 
$\lambda_{a_ib_i}$ to $\rho_{a_ib_i}$, whereas $\beta_i$ commutes with $\lambda_{a_ib_i}$.
Since the images of $\e_{a_i}\e_{b_i}$ and $\beta_i$ under $\phi$ are the same, this gives
\begin{equation*}
\phi(\rho_{b_ia_i}) = \phi(\lambda_{b_ia_i}).
\end{equation*}
As in the proof of Proposition \ref{p:quots}, we appeal to \cite{Ger84} to deduce that
$\phi$ factors through $\Sautn\to{\rm{SL}}(n,\Z)$.   

\medskip
 Next we consider the effect of the relations (3) and (4) from Lemma \ref{3lem}.
Unfortunately, these
are relations in \autn\ not \sautn. But 
since $R_i$ commutes with $\e_{a_j}$ when $j\neq i$ we have the following
relation in \sautn,
\begin{equation*}
R_i^{-1}\e_{a_i}\e_{a_j} R_i\e_{b_i}\e_{a_j}=\lambda_{a_ib_i}^2.
\end{equation*}
If $t=R_i$ then applying $\phi$ to this equation gives  $\phi(e_{a_i}e_{b_i})=\phi(\lambda_{a_ib_i})^2$.
Conjugating both sides by the permutation $(a_i\ a_j)(b_i\ b_j)$, we get the same
equality with $j$ subscripts. Since all the automorphisms with $i$ subscripts commmute with
those that have $j$ subscripts, we deduce 
$$
\hskip -1.7in (*)\hskip 2in  \phi( \e_{b_i}\e_{b_j} \e_{a_i}\e_{a_j})
= \phi(\lambda_{a_ib_i}^2\lambda_{a_jb_j}^2).
$$

If $t=R_iR_jv$ for some $j\neq i$ and $v$ a (possibly empty) word in the $R_k$ with $k\neq i,j$, then combining relation (4) for $i$ and $j$ gives
$$t^{-1}e_{a_i}e_{a_j}te_{b_i}e_{b_j}= R_i^{-1}R_j^{-1}e_{a_i}e_{a_j}R_jR_ie_{b_i}e_{b_j}= \lambda_{a_ib_i}^2 \lambda_{a_jb_j}^2.$$
Applying $\phi$ to this equation gives equation $(*)$ in this case as well.

If $t=R_iR_j^{-1}v$, then relation (3) for $i$ and relation (4) for $j$ give 
$$te_{b_i}e_{a_j}t^{-1}e_{a_i}e_{b_j}= R_iR_j^{-1}e_{b_i}e_{a_j}R_i^{-1}R_je_{a_i}e_{b_j}= \lambda_{a_ib_i}^2 \lambda_{b_ja_j}^2.$$
Applying $\phi$ to this equation gives 
$
\phi( \e_{a_i}\e_{b_j} \e_{b_i}\e_{a_j})
= \phi(\lambda_{a_ib_i}^2\lambda_{b_ja_j}^2).
$
Conjugating both sides by   $(a_j\, b_j)\e_{a_k}$ for some $k\neq i,j$ gives equation $(*)$ once again.

\medskip
Next we claim that equation $(*)$ forces $\phi$ to factor not only through 
$\Sautn\to {\rm{SL}}(n,\Z)$ but also through 
$\Sautn\to {\rm{SL}}(n,\Z_2)$. In order to prove this, it suffices to argue that
the image under $\phi$ of some Nielsen transformation has order at most 2.

Let $\overline\alpha$ denote the image of $\alpha\in{\rm{SAut}}(F_n)$ in
\sln. Consider the subgroup ${\rm{SL}}(4,\Z)\subset {\rm{SL}}(n,\Z)$ 
corresponding to the sub-basis $\{a_i,b_i,a_j,b_j\}$.
Equation $(*)$ tells us that $\overline\lambda_{a_ib_i}^2\overline\lambda_{a_jb_j}^2$
becomes central
in the image of ${\rm{SL}}(4,\Z)$ under $\phi$.
But in this copy of ${\rm{SL}}(4,\Z)$
one has the relations 
$[\overline\lambda_{a_ib_i}^2,\overline\lambda_{b_ib_j}]=
\overline\lambda_{a_ib_j}^2$ and $[\overline\lambda_{a_jb_j}^2,\overline\lambda_{b_ib_j}]=1$.
So forcing $\overline\lambda_{a_ib_i}^2\overline\lambda_{a_jb_j}^2$ to become central 
implies that $\phi(\lambda_{a_ib_j})^2=1$, as required.
 
\medskip
We have proved that $\phi$ factors through 
$\Sautn\to {\rm{SL}}(n,\Z_2)$. The final point 
 to observe is that the restriction to $T$ of this last map is injective; in
particular the image of $t$ is non-trivial, and hence so is the image of 
$\ker \phi$. Thus the image of $\phi$ in $G$ is a 
proper quotient of the {\em{simple}} group ${\rm{SL}}(n,\Z_2)$, and therefore is trivial.
\end{proof}

\section{Actions on generalized spheres and acyclic homology manifolds}\label{s:sphere}

Because the fixed point set of a finite-period homeomorphism of a sphere or contractible manifold need not be a manifold, we must expand the category we are working in to that of generalized manifolds.  We follow the exposition in Bredon's book on Sheaf Theory \cite{BredonBook}.  All homology groups in this section are Borel-Moore homology with compact supports and coefficients in a sheaf $\A$ of modules over a principle ideal domain $L$.  The homology groups of $X$ are denoted $H^c_*(X;\A)$. 
If $X$ is a locally finite CW-complex and $\A$ is the constant sheaf $X\times L$ (which we will denote simply by $L$), then 
$H^c_*(X;L)$ is isomorphic to singular homology with coefficients in $L$ (see 
\cite{BredonBook}, p.279).  

All cohomology groups are sheaf cohomology with compact supports, denoted $H_c^*(X;\A)$.  If $\A$ is the constant sheaf, this is isomorphic to \v Cech cohomology with compact supports.  If $F$ is a closed subset of $X$, then sheaf cohomology satisfies $H^k_c(X,F;\A)\cong H^k_c(X\smallsetminus F;\A)$.

In fact, the only sheaves we will consider other than the constant sheaf are the   sheaves $\mathcal O_k$  associated to the pre-sheaves $U\mapsto  H_k^c(X,X\smallsetminus U;L)$.  

\subsection{Homology manifolds} 
Let $L$ be one of $\Z$ or $\Z_p$ (the integers mod $p$,
where $p$ is a prime).

\begin{definition} \label{d:gm}(\cite{BredonBook}, p.329) 
 An {\em{$m$-dimensional
homology manifold over $L$}} (denoted $m$-{\rm{hm}}$_L$) is
a locally compact Hausdorff space  $X$  with finite homological dimension over  $L$, that has the local homology properties of a manifold.  Specifically, the sheaves $\Oh_k$ are locally constant with stalk $0$  if $k\neq m$  and $L$ if $k=m$. The sheaf $\Oh=\Oh_m$ is called the {\em orientation sheaf}.  
\end{definition} 
 

We will further assume that our homology manifolds are first-countable.  

\begin{definition}
If $X$ is an $m$-{\rm{hm}}$_L$ and $H^c_*(X; L) \cong H^c_*( \bS^m; L)$ then $X$
 is called a {\em{generalized $m$-sphere}} over $L$.
 \end{definition}
 
 \begin{definition} If $X$ is an $m$-{\rm{hm}}$_L$  with $H^c_0(X;L)=L$ and $H_k^c(X;L)=0$ for  $k>0$, then $X$ is said to be {\em $L$-acyclic}.
 \end{definition}

There is a similar notion of {\em cohomology manifold over $L$}, denoted $m$-{\rm{cm}}$_L$ (see \cite{BredonBook}, p. 373).  If $L=\Z_p$,  a connected space $X$ is an $n$-{\rm{cm}}$_{L}$ if and only if it is an $n$-{\rm{hm}}$_{L}$ and is  locally connected (\cite{BredonBook}, p.~375 Theorem 16.8 and footnote).  If $X$ is a locally connected homology manifold over $\Z_p$, then the fixed point set of any homeomorphism of order $p$ is also locally connected (see \cite{Borel}, Theorem 1.6, p. 72, where there is a stronger connectivity statement ($clc_L$), but the proof, which relies on Prop. 1.4, p. 68, also applies to local connectivity).  These remarks show that the theorems we state below for  homology manifolds  are also valid for cohomology manifolds.

 Finally, we note that homology manifolds satisfy Poincar\'e duality between Borel-Moore homology and sheaf cohomology (\cite{BredonBook}, Thm 9.2), i.e.
if $X$ is an $m$-hm$_L$ then 
$$H^c_k(X;L)\cong H^{m-k}_c(X;\Oh).$$

\subsection{Elements of Smith Theory} 

There are two types of Smith theorems, usually referred to as ``global" and ``local" Smith theorems.   The global theorems require only that $X$ be a locally compact Hausdorff space with the homology of a sphere or a point, while the local theorems concern homology manifolds.  These were originally proved by P.~A.~Smith  (\cite{SmithI},\cite{SmithII}), but we follow the exposition in  Bredon's book and Borel's Seminar on Transformation groups \cite{Borel}.

\begin{theorem}[The Local Smith Theorem, \cite{BredonBook}, Thm 20.1, Prop 20.2, pp. 409-410] Let $p$ be a prime and $L=\Z_p$.
The fixed point set of any action of $\Z_p$ on an $n$-hm$_L$  is the disjoint union of (open and closed) components each of which is
an $r$-{\rm{hm}}$_L$ with $r\leq m$. If  $p$ is odd then each component of the 
fixed point set has even codimension.
\end{theorem}

By invariance of domain  for homology manifolds (\cite{BredonBook}, Corollary 16.19, p. 383) the fixed point set of any {\it non-trivial} action of $\Z_p$ on a connected, locally connected $m$-{\rm{hm}}$_{\Z_p}$ 
 is  a (locally connected) $r$-{\rm{hm}}$_{\Z_p}$ with $r \leq m-1$.

 \begin{theorem}[Global Smith Theorems, \cite{BredonBook}, Corollaries 19.8 and 19.9, p. 144]\label{GST} Let $p$ be a prime and $X$ a locally compact Hausdorff space of finite dimension over $\Z_p$.  Suppose that $\Z_p$ acts on $X$ with fixed point set $F$.
 \begin{itemize}
 \item If $H^c_*(X;\Z_p)\cong H^c_*(\bS^m;\Z_p)$, then  $H^c_*(F;\Z_p)\cong H^c_*(\bS^r;\Z_p)$ for some $r$ with $-1\leq r\leq m$. If $p$ is odd, then $m-r$ is even.
 \item If $X$ is $\Z_p$-acyclic, then $F$ is $\Z_p$-acyclic (in particular non-empty and connected).
\end{itemize}
\end{theorem}

In section 19 of \cite{BredonBook} the Global Smith Theorem is stated
for cohomology; the homology version above follows using the Smith
Theory sequence (132) on page 408 of \cite{BredonBook}.  The details
of this translation have been worked out by Olga Varghese in
\cite{Varghese}.

Together these theorems imply 
 
\begin{corollary}\label{smith}Let $X$ be an $m$-{\rm{hm}}$_{\Z_p}$.
\begin{itemize}
\item If $X$ is a generalized $m$-sphere over $\Z_p$, the fixed point set of any homeomorphism of order $p$ is a (possibly empty) generalized $r$-sphere, 
with $r\leq m-1$.  If $p$ is odd, $r\leq m-2$.
\item If $X$ is  $\Z_p$-acyclic, the fixed point set of any homeomorphism of order $p$ is a (non-empty) $\Z_p$-acyclic $r$-{\rm{hm}}$_{\Z_p}$, for some $r\leq m-1$.  If $p$ is odd, $r\leq m-2$.
\end{itemize}
\end{corollary}

We want to use this corollary as the basis for an
induction that bounds the dimensions in which elementary
$p$-groups can act effectively on generalized spheres
and acyclic homology manifolds. But in the case of 
spheres we need an additional result that guarantees the existence
of fixed points. This is provided by P.A.~Smith's theorem that 
$\Z_p\times\Z_p$ cannot act freely on a generalized sphere over $\Z_p$
(see \cite{SmithPNAS}; cf
Theorem \ref{fthree} below).
 
The proof of the following theorem is again
due to P.A.~Smith \cite{SmithPNAS}. (In \cite{SmithPNAS} he only gave the
proof for generalized spheres, but the 
acyclic case is similar.)  

\begin{thm}\label{Smith}  
If $m<d-1$, the group $(\Z_2)^d$ cannot act effectively
on a generalized $m$-sphere over $\Z_2$ or a $\Z_2$-acyclic $(m+1)$-dimensional homology manifold over $\Z_2$.

If $m<2d-1$ and $p$ is odd,  then $(\Z_p)^d$ cannot act effectively
a generalized $m$-sphere or a $\Z_p$-acyclic $(m+1)$-dimensional homology manifold over $\Z_p$.
\end{thm}

\begin{proof} The cases that arise when $d=1$ are vacuous or trivial
except when $p$ is odd and the putative action is on a $1$-{\rm{hm}}$_{\Z_p}$,
in which case
 one needs to recall that a $1$-{\rm{hm}}$_{\Z_p}$ is an actual
manifold.

We assume $d\ge 2$ and proceed by induction.
Let $X$ be one of the spaces that the theorem asserts 
$G:=(\Z_p)^d$ cannot act effectively on. 

Among the  non-trivial elements of $G$ we choose one, $a$ say,
whose  fixed point
set $F_a$ is maximal with respect to inclusion. We also
choose a complement $G_0\cong(\Z_p)^{d-1}$ to $\langle a\rangle$
in $G$. From Thorem \ref{GST} (in the acyclic case) and Smith's theorem
for $\Z_p\times\Z_p$ (in the case of spheres), we know that $F_a$
is non-empty.
We shall prove that
$F_a=X$ by assuming it false
and obtaining a contradiction. 

If 
$F_a$ is not the whole of $X$ then it
has codimension at least 1 if $p=2$ and codimension at least 2 if $p$ is odd. In the light of 
Corollary \ref{smith}, we may apply induction to the action of $G_0$
on $F_a$ and hence conclude that some non-trivial element $b\in G_0$
fixes $F_a$ pointwise; in other words $F_a\subseteq F_b$. But $F_a$
is maximal, so   $F_b=F_a$, which implies $F_a=Fix(A)$ for $A=\langle  a,b\rangle$.
Thus for any non-trivial element  $x$ of $A$
we have  $F_a\subseteq F_x$, so again maximality tells us
that $F_a=F_x$.  

Theorem 4.3 on page 182 of \cite{Borel}  (which requires us to know
that $\fix(A)$ is non-empty) provides a formula
relating the dimensions of the fixed point sets of elements of $A$:
writing $n={\rm{dim}}_p(X), r={\rm{dim}}_p(\fix(A))$, and $ r_C={\rm{dim}}_p(\fix(C))$ for each cyclic subgroup $C< A$, we have
$$n-r=\sum (r_C-r),$$
where the sum is taken over the non-trivial cyclic subgroups of $A$. 
We have just argued that $F_a=\fix(C)=\fix(A)$ for all non-trivial $C<A$,
so  each summand on the right is 0 and hence $n=r={\rm{dim}}_p(F_a)$.
Since $X$ is connected, invariance of domain gives $X=F_a$, i.e. $a$ acts trivially. This contradiction completes the induction.
\end{proof}

We need one more result from Smith theory:

\begin{theorem}\label{fthree} Let $X$ be a generalized sphere over $\Z_2$  or a $\Z_2$-acyclic {\rm{hm}}$_{\Z_2}$, and let
$a$ and $b$ be commuting 
 homeomorphisms of $X$, each of order $2$, with fixed point sets
 $F_a$ and $F_b$. If $F_a=F_b$ then $a=b$.
 \end{theorem}

\begin{proof} For actions on generalized spheres, this is explicit in
\cite{SmithPNAS}, so we consider only the acyclic case.

If $a\neq b$ then the subgroup $A\leq 
{\rm{Homeo}}(X)$ generated by $a$ and $b$ is isomorphic to $\Z_2\times \Z_2$ and $\fix(a)=\fix(b)=\fix(A)$. Thus in the formula
 $n-r=\sum (r_C-r)$ displayed in the preceding proof, the only non-zero
 summand on the right is the one for $\langle ab\rangle$. Hence $n=r_C$,
 that is,
 ${\rm{dim}}_p(X)={\rm{dim}}_p(\fix(ab))$.  Since $X$ is connected, invariance of domain gives $X=\fix(ab)$, which means that $a=b$.  
\end{proof}

\subsection{Actions on generalized spheres and $\Z_3$-acyclic homology manifolds over $\Z_3$, for $n$ even}\label{ss:Zthree}

The results we have developed to this point easily yield the following theorem, for $n$ even.

\begin{theorem} Let $X$ be a generalized $m$-sphere over $\Z_3$
or a $\Z_3$-acyclic $(m+1)$-dimensional
homology manifold  over $\Z_3$, and let $\phi\colon\Sautn\to{\rm{Homeo}}(X)$ 
be an action.  If $n$ is even and $m<n-1$, then $\phi$ is trivial.  
\end{theorem}
\begin{proof}
Write $n=2d$ and
let $T\subset{\rm{SAut}}(F_{2d})$ be as in Lemma \ref{whatTis}.
Since  $T\cong(\Z_3)^d$ and $m<n-1=2d-1$,  Theorem \ref{Smith} tells us that 
 $T$ cannot act effectively on $X$, so
$\phi(t)=1$ for some $t\in T\smallsetminus\{1\}$. We proved
in Proposition \ref{p:3quots} that this forces $\phi$ to be the
trivial map. 
\end{proof}

\subsection{Actions on generalized spheres and $\Z_2$-acyclic homology manifolds over $\Z_2$} 

The proof of Theorems \ref{t:2sphere} and \ref{t:2acyclic} is considerably
more involved than that of the preceding result. This is largely due to the
fact that Corollary \ref{smith} yields a weaker conclusion for $p=2$ than
for odd primes. Lemma \ref{l:tech} will allow us to circumvent this difficulty.
It relies on the separation property of codimension 1 fixed point sets
that is established in Lemma \ref{complement} using
Poincar\'e duality and the following theorem about sheaf cohomology. 

\begin{theorem}
[Theorem 16.16, \cite{BredonBook}] If $X$ is a connected $m$-{\rm{cm}}$_L$ with orientation sheaf $\Oh$, and $F$ is a proper closed subset, then for any non-empty open subset $U$ \begin{enumerate}
\item $H^m_c(U;\mathcal O)$ is the free $L$-module on the components of $U$
\item $H^m_c(F;L)=0$
\end{enumerate}
\end{theorem}

\begin{lemma}\label{complement} Let  $X$ be a  
generalized $m$-sphere over $\Z_2$ or a
$\Z_2$-acyclic $m$-{\rm{hm}}$_{\Z_2}$, and let $\tau$ be an involution of $X$.
 If $\fix(\tau)$ has dimension $m-1$, then $X\smallsetminus\fix(\tau)$ 
has two $\Z_2$-acyclic components and $\tau$ interchanges them. 
\end{lemma}
\begin{proof} If $m=1$, then $X$ is a circle or a line, $\fix(\tau)$ is two points or one, and the theorem is clear, so we may assume $m\geq 2$.  Let $F=\fix(\tau)$, and set $L=\Z_2$.  Since $F$ is closed, the long exact sequence in sheaf cohomology for the pair $(X,F)$ reads
$$\ldots\to \Hc^{m-2}(F;L)\to \Hc^{m-1}(X\smallsetminus F;L)\to \Hc^{m-1}(X;L)\to \Hc^{m-1}(F;L)\to $$
$$H^{m}_c(X\smallsetminus F;L)\to \Hc^{m}(X;L)\to\Hc^{m}(F;L)\to 0$$
By (2) above, the last term $\Hc^{m}(F;L)$  is $0$.  Since $L=\Z_2$, the orientation sheaf $\Oh$ is actually constant, and by (1), we get $\Hc^{m}(X;L)=\Hc^{m-1}(F;L) = L$.  

Poincar\'e duality says $H^{k}_c(X;L)\cong H_{m-k}^c(X;L)$;  in particular, $H^{m-1}_c(X;L)\cong H_1^c(X;L)=0$ (since $m\geq 2$), and the end of the sequence is
$$ 0\to L\to \Hc^m(X\smallsetminus F;L)\to L\to 0$$
Thus  $\Hc^m(X\smallsetminus F;L)\cong L\oplus L$, and another application of (1) shows that $X\smallsetminus F$ has two components.  (This is the argument in \cite{BredonBook}, Cor. 16.26.)

Suppose $X$ is $L$-acyclic. Applying Poincar\'e duality to each remaining term in the long exact sequence ($F$ has dimension $m-1$) gives
$$\ldots \to \Hcp_{k}(F;L)\to \Hcp_{k}(X\smallsetminus F;L)\to \Hcp_{k}(X;L)\to \ldots$$
for $k\geq 1$.  Since $F$ and $X$ are $L$-acyclic, this shows that each component of $X\smallsetminus F$ is also $L$-acyclic.

If $X$ is a generalized $m$-sphere then $F$ is a generalized $(m-1)$-sphere, and the above argument shows that most of the homology of $X\smallsetminus F$ vanishes as in the acyclic case. In dimensions $m$ and $m-1$ we have 
$$0\to \Hcp_{m}(X\smallsetminus F;L)\to \Hcp_{m}(X;L)\to\Hcp_{m-1}(F;L)\to \Hcp_{m-1}(X\smallsetminus F;L)\to 0$$
which becomes 
$$0\to \Hcp_{m}(X\smallsetminus F;L)\to L\cong   L\to \Hcp_{m-1}(X\smallsetminus F;L)\to 0$$ so again the homology of $X\smallsetminus F$ vanishes in positive degrees, and each component of $X\smallsetminus F$ is acyclic. 

In both situations the complement of $F$ has two $\Z_2$-acyclic components.  Since the involution acts freely on this complement, it cannot preserve either component, by the Global Smith Theorem. 
\end{proof}
 
\begin{lemma}\label{l:tech} Let $X$ be a generalized $m$-sphere 
over $\Z_2$  or  a $\Z_2$-acyclic $(m+1)$-{\rm{hm}}$_{\Z_2}$, and let $G$ be a group acting by homeomorphisms on $X$.   Suppose $G$ contains a subgroup $P\cong \Z_2\times\Z_2$ all of whose non-trivial elements are conjugate in $G$.  If 
 $P$ acts non-trivially, then the fixed point sets of 
 its non-trivial elements have  codimension at least $2$, and  $m\ge 2$.
 \end{lemma}
 
 \begin{proof}  Since the non-trivial elements of $P$ are all conjugate, they must all act non-trivially.
 
 Let $a$ and $b$ be generators of $P$.  If
 $\fix(a)$ had codimension 1, then by Lemma~\ref{complement}, its complement in $X$ would have two components
 and the action of $a$ would interchange these. Consider the action of 
 $b$: since it commutes with $a$ it leaves $\fix(a)$
 invariant, so it either interchanges the components of the complement
 or leaves them invariant. Reversing the roles of $b$ and
 $ab$ if necessary, we may assume that it interchanges them
 and hence that 
 $\fix(b)\subset\fix(a)$. Since $a$
 and $b$ are conjugate, invariance of domain for 
 homology manifolds implies that $\fix(a)=\fix(b)$
 and hence, by Theorem \ref{fthree}, that the  actions of $a$ and $b$
 on $X$ are identical. Thus $ab$ acts trivially,  contradicting the assumption that the action of $P$ is non-trivial.

Thus the fixed point set of any non-trivial element of $P$ has codimension at least 2.    If $m=1$ and $X$ is a generalized sphere this says $\fix(a)=\fix(b)=\emptyset$.  If $X$ is 2-dimensional and acyclic,  then $\fix(a)$ and $\fix(b)$ are $0$-dimensional acyclic homology manifolds, i.e. points, so $\fix(a)\subset \fix(b)$ implies $\fix(a)=\fix(b)$.  In either case,  Theorem $\ref{fthree}$ again implies that $ab$ acts trivially, contradicting our assumptions.
 \end{proof}

\begin{prop}\label{p:sln} Let $X$ be a generalized $m$-sphere or  
$\Z_2$-acyclic $(m+1)$-{\rm{hm}}$_{\Z_2}$. 
If $m<n-1$ and $n\ge 3$, then any action of 
${\rm{SL}}(n,\Z_2)$ on $X$ is trivial.
\end{prop}

\begin{proof} 
Since ${\rm{SL}}(n,\Z_2)$ is simple, it is enough to find a
subgroup of ${\rm{SL}}(n,\Z_2)$ that cannot act effectively on
$X$.

The elementary matrices $E_{j1},\ j\neq 1,$ generate an elementary
2-group $Q\cong (\Z_2)^{n-1}$. All elementary matrices are conjugate
in ${\rm{SL}}(n,\Z_2)$. Moreover $E_{32}E_{21}E_{32}=E_{31}E_{21}$.
Thus we are in the situation of Lemma \ref{l:tech} with $a=E_{21}$
and $b=E_{31}$. An appeal to that lemma completes the proof in the
case $n=3$.

If $n\ge 4$ then  ${\rm{SL}}(n,\Z_2)$ contains a larger elementary
2-group than $Q$, namely that generated by the elementary matrices
$E_{ij}$ with $i\le n/2$ and $j> n/2$. This has rank at least $n$, so
Theorem \ref{Smith} tells us it cannot act effectively $X$.
\end{proof}

Proposition~\ref{p:quots} and Proposition~\ref{p:sln} together give:

\begin{corollary}\label{c:saut}  Let $X$ be a generalized $m$-sphere over
$\Z_2$ or  a
$\Z_2$-acyclic $(m+1)$-{\rm{hm}}$_{\Z_2}$, with $m<n-1$.  If a non-central element of $W_n$ is in the kernel of an action of \sautn\ on $X$, then the action is trivial. 
\end{corollary}

\subsection{Proof of Theorems \ref{t:2sphere} and \ref{t:2acyclic}}\label{p:2sphere}
 
We retain the notation introduced at the beginning of section \ref{s:quots}.
 
 Let $X$ be a generalized $m$-sphere or  a $\Z_2$-acyclic $(m+1)$-dimensional homology manifold  over $\Z_2$, with $m<n-1$.  Let $\Phi:{\rm{SAut}}(F_n)\to \rm{Homeo}(X)$ be an action of
\sautn\  on $X$.
In the light of the preceding corollary, we will be done if we can prove that
the kernel of $\Phi$ contains an element of  
$\SN\cong (\Z_2)^{n-1}$ other than $\Delta=e_1\dots e_n$.
 
\medskip
 
If $n=3$, then   
conjugating
$a:=e_1e_2$ by $(1\ 3)e_2$ and $(2\ 3)e_1$, respectively, yields
$b:=e_2e_3$
and $ab = e_1e_3$. Thus we may appeal to Lemma \ref{l:tech},
to see that the action of $\SN$ on $X$ is trivial if $m<2$. 

\smallskip

If $n=4$, then $\SN$ is generated by $a,\ b,$ and $c:=e_2e_4$,
which are conjugate in $\SW$ to each other and to each of the
products $ab, ac$ and $bc$. If the action of $\SN$ on $X$ is
trivial then we are done. Suppose that this is not
the case. 
We know from
 Lemma~\ref{l:tech} that $\fix(a)$ is a generalized $d$-sphere over 
 $\Z_2$ or a $\Z_2$-acyclic $(d+1)$-{\rm{hm}}$_{\Z_2}$
 with $d<1$. Since $b$ and $c$ commute with $a$, the group $\langle b,c\rangle\cong \Z_2^2$ acts on $\fix(a)$, so by Theorem~\ref{Smith} some element acts trivially, say $g$.  Then $\fix(g)\supset \fix(a)$. But since $a$ and $g$ are conjugate in $\SW$, this implies $\fix(a)=\fix(g)$, and then by Theorem~\ref{fthree}, $ag$ acts trivially on $X$.  
If $g=b$ or $g=c$, we have found a non-central element of the kernel of $\Phi|_{SW_n}$, so $\Phi$ is trivial by Corollary~\ref{c:saut} .  If $g=bc$, 
the action factors through ${\rm{PSL}}(4,\Z)$, which contains a subgroup isomorphic to $(\Z_2)^4$ generated by $e_1e_2,\ e_2e_3,\ \sigma=(12)(34),$ and $\tau=(13)(24)$.  By Theorem~\ref{Smith}, some nontrivial element of this subgroup must map trivially to $\rm{Homeo}(X)$.  Pulling back these elements to $\SW$, we see that 
some element of the form $e_ie_j$ or $\gamma$ or  $e_ie_j\gamma$, with $\gamma
\in\langle\sigma,\tau\rangle$,
 is in the kernel.   Corollary~\ref{c:saut} again shows that $\Phi$ is trivial.
 
 \medskip
 
Now we suppose $n>4$ and proceed by induction.
If $e_1e_2$ acts trivially then we are done by Corollary~\ref{c:saut}.
If not then, appealing to Lemma \ref{l:tech} once more,
we may suppose that the fixed point set of $e_1e_2$ in $X$ is a 
generalized
$r$-sphere or $\Z_2$-acyclic $(r+1)$-homology manifold
over $\Z_2$  with $r<n-3$; call it $Y$.
The centralizer of $e_1e_2$ will
act on $Y$. This centralizer contains a  copy of ${\rm{SAut}}(F_{n-2})$
corresponding to the sub-basis $a_3,\dots,a_n$, 
and by induction this acts trivially on $Y$. In particular,
 the fixed point set of $e_3e_4$ contains that of $e_1e_2$. Similarly,
the reverse inclusion holds. But then by Theorem~\ref{fthree} 
the actions of $e_3e_4$ and $e_1e_2$ on  $X$ must be the same.
Thus the kernel of any homomorphism ${\rm{SAut}}(F_n)\to \rm{Homeo}(X)$ intersects
$N$ in $e_1e_2e_3e_4\neq\Delta$, and   Corollary \ref{c:saut} says that
the action is therefore trivial. \qed

\begin{remark} For $n=3$,  Theorem \ref{t:2sphere}
 also follows from the results of \cite{BriVog03:2} because a generalized $\Z_2$-sphere
of dimension one is just a circle and it is shown in \cite{BriVog03:2}  that any action of 
\sautn\ by homeomorphisms on a circle is trivial for $n\geq 3$.
\end{remark}

\begin{remark}\label{ParFlaw} This work was stimulated in part
by the proof of Corollary \ref{c:sl} suggested by 
Zimmermann in  \cite{Zim06}. His proof relied on earlier work of
Parwani \cite{Par06} which sets forth a good strategy but contains a 
flaw: it is assumed in \cite{Par06}
that if $X$ is a 
homology manifold over $\Z$ with the $\Z_2$-homology of a sphere,
then the fixed point set of any involution of $X$ will again be such a space; this is false (see the following remark).
It is also assumed in \cite{Par06} that such a fixed point
set will be an ENR, and this is also false.
\end{remark}

\begin{remark} In \cite{Jones}, L.~Jones showed that almost any PL homology manifold over $\Z_2$ satisfying the Smith conditions can arise as the fixed point set of  
a involution of a genuine sphere.   In particular, the fixed point set of an involution of a sphere need not be a $\Z$-homology manifold.    

There are also involutions of spheres for which the fixed point set is not locally 1-connected, so in particular is not an
ENR.  Indeed Ancel and Guilbault \cite{AG}
proved that if $\overline M = M\cup \Sigma$ is {\em{any}} $Z$-set compactification of
a contractible $n$-manifold $M$, with $n>4$, then the double of $\overline M$ along $\Sigma$ is homeomorphic to
the $n$-sphere. One can  realize $\Sigma$ as the fixed point set of the involution that interchanges the
two copies of $M$ in this double, and $\Sigma$ need not be locally 1-connected. To obtain a concrete example, we can
take $M$ to be the universal cover of one of the aspherical manifolds constructed by Davis \cite{davis}
and  take $\Sigma$ to be its ideal boundary (cf.~\cite{fish}).
\end{remark}
    
\section{Actions on arbitrary compact homology manifolds}

In \cite{MannSu}, Mann and Su use Smith theory and a spectral sequence
developed by Swan \cite{swan} to prove that for every prime $p$ and
every compact $d$-dimensional homology manifold $X$ over $\Z_p$,
the sum of whose mod $p$ betti numbers is $B$, there
exists an integer $\nu(d,B)$, depending only on $d$ and $B$,
so that  
$\Z_p^r$ cannot act effectively by homeomorphisms
on $X$  if $r>\nu(d,B)$. (An explicit bound on $\nu$ is given.)

If $n$ is sufficiently large then the
 alternating group $A_n$ will be simple and contain a copy
 of $\Z_p^\nu$. Hence it will admit no non-trivial action on $X$.
Theorem \ref{noAct} stated in the introduction is an immediate consequence of this result and 
Proposition \ref{p:quots}, since $SW_n\subset {\rm{SAut}}(F_n)$ 
contains a copy of $A_n$. \hfill $\square$

\smallskip

The preceding argument allows one to bound the constant
$\eta(p,d,B)$ in Theorem \ref{noAct} by a multiple (depending
on $p$) of $\nu(d,B)$. In
the cases $p=2$ and $p=3$ one can sharpen this estimate by 
appealing directly to Propositions \ref{p:quots}
and   \ref{p:3quots} instead of using $A_n$.

\begin{remark}
Various of the Higman-Thompson groups, including
Richard Thompson's {\em{vagabond group}} $V$, are finitely presented,
simple, and contain  an isomorphic copy of every finite group \cite{higman}.
Given any class of objects each of which has
the property that some finite group cannot act effectively on it, groups
such as $V$ cannot act non-trivially on any object in the class.
In particular, it follows from the Mann and Su result that $V$ cannot
act non-trivially by homeomorphisms on any compact manifold.
And Theorem \ref{Smith} above implies that $V$ cannot act non-trivially
by homeomorphisms  on any finite-dimensional $\Z_p$-acyclic
homology manifold over $\Z_p$ for any prime $p$ (cf.~\cite{gang6}
and \cite{fisher}).
\end{remark}

\end{document}